\def\editmode{0}

\def\hyperlinks{1}  

\def\spsformat{0}

\def\singlenarrowcol{0} 

\def\bibfilenames{bibman_refs}

\def\journal{0}

\def\intermediatesteps{1}

\if\spsformat1

\documentclass{article}
\usepackage{spconf}

\else

\if\singlenarrowcol0
\documentclass[10pt,conference]{IEEEtran}

\IEEEoverridecommandlockouts
\else

\documentclass[onecolumn]{IEEEtran}
    \usepackage[pass]{geometry}
    \newgeometry{textwidth=252pt, textheight=672pt}

\fi
\fi

\usepackage{include_v8}
\usepackage{notation}

\newcommand{\conference}[1]{} 

\newcommand{\nextversion}[1]{}
\renewcommand{\define}{\triangleq}

\begin{document}

\title{Spatial Transformers for Radio Map Estimation}

\author{
    Pham Q. Viet and Daniel Romero
    \\
    Dept. of Information and Communication Technology,
    University of Agder, Grimstad, Norway.\\
    Email:\{viet.q.pham,daniel.romero\}@uia.no.
}
\maketitle

\begin{abstract}
    Radio map estimation (RME) involves spatial interpolation of radio measurements to predict metrics such as the received signal strength at locations where no measurements were collected. The most popular estimators nowadays project the measurement locations onto a regular grid and complete the resulting measurement tensor with a  convolutional deep neural network. Unfortunately, these approaches suffer from poor spatial resolution and require a very large number of parameters. The first contribution of this paper addresses these limitations by means of an attention-based  estimator named Spatial TransfOrmer for Radio Map estimation (STORM). This scheme not only outperforms the existing estimators, but also exhibits lower computational complexity, translation equivariance, rotation equivariance, and full spatial resolution. The second contribution is an extended transformer architecture that allows STORM to perform active sensing, by which the next measurement location is selected based on the previous measurements. This is particularly useful for minimization of drive tests (MDT) in cellular networks, where operators request  user equipment to collect measurements. Finally, STORM is extensively validated by experiments with one ray-tracing and two real-measurement datasets.

\end{abstract}

\newcommand\blfootnote[1]{%
    \begingroup
    \renewcommand\thefootnote{}\footnote{#1}%
    \addtocounter{footnote}{-1}%
    \endgroup
}

\begin{IEEEkeywords}
    Radio map estimation, transformers, attention-based learning, deep learning, wireless communications. \blfootnote{    This work has been funded by the IKTPLUSS grant 311994 of the Research
        Council of Norway.}

\end{IEEEkeywords}

\begin{bullets}%

    \section{Introduction}
    \label{sec:intro}

    \blt[Overview]
    \begin{bullets}%
        \blt[Radio maps] Radio maps (see Fig.~\ref{fig:estimationexample}), also known as radio environment maps, provide radio frequency (RF) metrics such as the received signal strength across a geographical region~\cite{romero2022cartography,romero2024theoretical}.
        \blt[Applications]Radio maps find a large number of applications, including network planning, frequency planning, cellular communications, device-to-device communications, dynamic spectrum access,   robot path planning, aerial traffic management in unmanned aerial systems, and fingerprinting localization to name a few; see e.g.~\cite{\nextversion{abouzeid2013predictive,subramani2011practical,}cai2011ran,zalonis2012femtocell,romero2022cartography,teganya2019locationfree} and references therein.
    \end{bullets}%

    \blt[Related work]
    \begin{bullets}%
        \blt[Non-attention based estimators]
        \blt\cmt{RME SOA}
        \begin{bullets}%
            \blt[RME def]Radio map estimation (RME) involves constructing a radio map by relying on   measurements collected across the area of interest.
            \blt[estimators]
            \begin{bullets}%
                \blt[before DL]
                \begin{bullets}%
                    \blt[kernel]Before the advent of deep learning, the most popular estimators were built upon kernel-based learning (see~\cite{romero2017spectrummaps} and references therein),
                    \blt[kriging] Kriging~\cite{alayafeki2008cartography,\nextversion{agarwal2018spectrum,}shrestha2022surveying},
                    \blt[sparsity]sparsity-based inference~\cite{bazerque2010sparsity\nextversion{,bazerque2011splines,jayawickrama2013compressive}},
                    \blt[Matrix completion]matrix completion~\cite{schaufele2019tensor\nextversion{,khalfi2018airmap}},
                    \blt[dictionary learning] dictionary learning~\cite{kim2013dictionary}, and graphical models~\cite{ha2024location}.
                \end{bullets}%
                \blt[deep learning]The most recent estimators are based on   deep neural networks (DNNs); see e.g.~\cite{krijestorac2020deeplearning,levie2019radiounet\nextversion{,han2020power},teganya2020rme,shrestha2022surveying}. Unfortunately, these schemes entail grid discretization and a very large number of trainable parameters, which render them computationally expensive and drastically limit their spatial resolution. Besides, they lack important desirable properties in RME, such as translation and rotation equivariance.
            \end{bullets}%
        \end{bullets}%

        \blt[Active sensing] A relevant task in RME is also active sensing, where the next measurement location is to be decided given the previously collected measurements. The special case where the upcoming measurements must lie on the trajectory of a mobile robot such as a UAV has been considered in \cite{shrestha2022surveying}. Other works have proposed extensions and improvements in different settings, but the approaches therein are reminiscent of the aforementioned estimation schemes.
        %

        \blt[Papers related to RME that use transformers]Some works related to RME have  considered attention-based estimators, which is the topic of this paper.
        \begin{bullets}%
            \blt[Transformers on images]For example, vision transformers have been used to accommodate side information, such as
            \begin{bullets}%
                \blt[building maps]building maps~\cite{tian2021transformer,zheng2024site}
                \blt[satellite image]and satellite images~\cite{yu2024transformer}.
            \end{bullets}%
            \blt[Transformers on radio measurements]In other works, transformers are fed with radio measurements.
            \begin{bullets}%
                \blt[pandey2021limited]This is the case of
                \cite{pandey2021limited}, which
                uses a transformer for predicting what a device would measure given the measurements collected by a device with different hardware characteristics,
                \blt[wang2023bert] and of \cite{wang2023bert}, where  a transformer  fills missing RSS features for fingerprinting-based localization.
            \end{bullets}%
        \end{bullets}%
    \end{bullets}%
    \blt[Research gap]
    \begin{bullets}%
        \blt[RME]Thus, transformers have been applied to problems that are related to RME or to modified versions of the RME problem where transformers are used to process images.
        However, the plain RME problem, where a radio map needs to be constructed by relying solely on radio measurements and their locations, requires spatial interpolation of radio measurements and this has never been tackled using transformers.

        \blt[Active sensing]

    \end{bullets}%

    \begin{figure}
        \centering
        \includegraphics[clip, trim=3.8cm 0.4cm .8cm 1.8cm, width=.75\columnwidth]{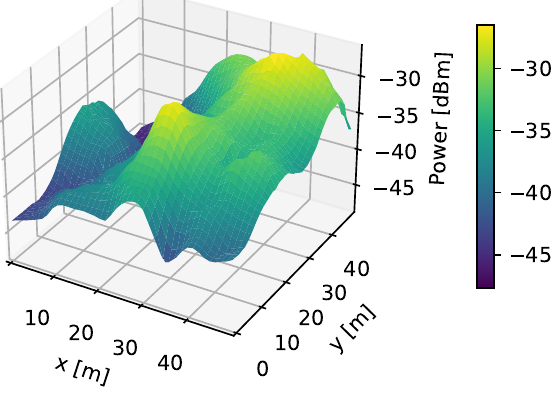}
        \caption{ Example of a radio map estimate obtained with the proposed STORM estimator in the USRP dataset; cf. Sec.~\ref{sec:simulations}. }
        \label{fig:estimationexample}
    \end{figure}

    \blt[Contributions]%
    \begin{bullets}%
        \blt[Attention-based RME]The first contribution of this paper is an attention-based  scheme referred to as \emph{Spatial TransfOrmer for Radio Map estimation} (STORM).
        \begin{bullets}%
            \blt[sets SoA performance]As shown by numerical experiments in one ray-tracing and two real-measurement datasets, STORM outperforms the existing estimators, thereby setting the state of the art in RME.
            \blt[relative to other estimators]
            \begin{bullets}%
                \blt[Relative to DNN estimators]Besides, it offers key advantages over existing DNN estimators:
                \begin{bullets}%
                    \blt[Complexity]First, its  complexity is significantly lower. For example, the RadioUnet estimators in  \cite{levie2019radiounet} have 6 M and 25 M parameters, whereas STORM has just 100 k parameters.
                    \blt[Gridless]Second,  STORM  is \emph{gridless}. This implies that, unlike previous DNN estimators, (i) its spatial resolution is not limited, (ii) it is translation and rotation equivariant, which are desirable properties in view of Maxwell's equations, and (iii), it need not be retrained e.g. if the grid spacing changes.
                    %
                    \blt[The map is estimated only where it is needed]Third, STORM can estimate the map only at the relevant locations, whereas previous DNN estimators need to compute the map at all grid locations.
                    \blt[Can accommodate measurements outside the region]Fourth, STORM can accommodate measurements outside the region of interest, as opposed to existing DNN estimators.
                \end{bullets}%
                \blt[Relative to Kriging]STROM also compares favorably to Kriging, which is the non-DNN estimator of choice
                \begin{journalonly}
                    , as empirically  corroborated in
                \end{journalonly}\cite{elfriakh2018crowdsourced,shrestha2024empiricaljpaper}. While Kriging requires a matrix inversion, whose complexity is cubic in the number of measurements, STORM enjoys quadratic complexity.

            \end{bullets}%

        \end{bullets}%
        \blt[Active sensing]The second contribution is an extension of STORM to active sensing. Given a set of candidate measurement locations, STORM will indicate which one of them should be selected to collect the next measurement so that the estimation error is approximately minimized.
        \begin{bullets}%
            \blt[Application: MDT]This is of special interest in the setup of \emph{minimization of drive tests} (MDT), a technology where cellular communication operators request  user equipment to collect geolocalized measurements.
            \blt[Transformer-based architecture]The approach relies on two interconnected encoder-decoder transformer networks.
        \end{bullets}%
    \end{bullets}%

    \blt[Paper structure] Sec.~\ref{sec:problem} formulates the RME and active sensing problems. Sec.~\ref{sec:transformers} reviews key notions about transformers. Secs.~\ref{sec:storm} and~\ref{sec:activestorm} propose STORM. Finally, Secs.~\ref{sec:simulations} and~\ref{sec:conclusions} respectively present the numerical experiments and conclusions.
    \blt[code]The code and data will be published on \url{www.radiomaps.org}.

    \blt[Notation]\emph{Notation}: Lowercase (uppercase) boldface letters represent column vectors (matrices). Equality by definition is denoted as $\define$. For a matrix $\bm A$, the $(m,n)$-th entry of $\SofMax(\bm A)$ is $\exp([\bm A]_{m,n})/\sum_{m}\exp([\bm A]_{m,n})$.

\section{The RME Problem}
\label{sec:problem}

This section presents the most prevalent formulation of  RME, both in the conventional  and active sensing setups. For simplicity,
the signal strength is quantified here by the
received signal power, but other metrics can readily be adopted.

\begin{bullets}
    \blt[model]
    \begin{bullets}
        \blt[region] Let $\area\subset\rfield^{\areadim}$ encompass the Cartesian
        coordinates of all points within the region of interest,
        whose dimension $\areadim$ is  typically  2 or 3. Very often, $\area$ is a
        rectangular area in a horizontal plane.
        \blt[map]A \emph{power map} is a function that returns the  signal power
        $\pow(\loc)$ that a sensor with an isotropic antenna at location
        $\loc \in \area$ would receive. This power is the result of the
        contribution of one or multiple transmitters as well as the propagation effects in the environment.
        %

        \blt[measurement model] The received power is measured at $\measnum$
        locations $\left\{\loc_{\measind} \right\}_{\measind=1}^{\measnum}
            \subset \area$ by one or
        multiple receivers (or sensors). The $\measind$-th measurement can be written as
        $
            \powmeas[\measind] = \pow \left( \loc_{\measind} \right) + \measnoise_{\measind}$,
        where $\measnoise_{\measind}$ denotes measurement error. 
    \end{bullets}
    \blt[Problem Formulations]
    \begin{bullets}%
        \blt[RME]
        \begin{bullets}%
            \blt[given] Given
            $\left\{(\loc_{\measind},
                \powmeas[\measind])\right\}_{\measind=1}^{\measnum}$,
            \blt[requested]the RME problem is to estimate  $\pow(\loc)$, $\loc\in\area$.
        \end{bullets}%
        \blt[Active Sensing]On the other hand, in the active sensing problem,
        \begin{bullets}%
            \blt[given]one is given the measurements
            $\left\{(\loc_{\measind},
                \powmeas[\measind])\right\}_{\measind=1}^{\measnum}$ as well as the set $\{\loc_{\measind}\}_{\measind=\measnum+1}^{\measnum+\CanNum}$ of candidate locations
            \blt[requested] and selects one of these candidate locations, say $\loc_{\CanInd}$, to collect the next measurement $\powmeas[\CanInd]$. The goal is to choose $\CanInd$ so that a target metric of the estimation error is minimized given the measurements $\{(\loc_{\measind},
                \powmeas[\measind])\}_{\measind=1}^{\measnum} \cup\{
                (\loc_{\CanInd}, \powmeas[\CanInd])
                \}$.

        \end{bullets}

    \end{bullets}
\end{bullets}

\section{Background on Transformers}
\label{sec:transformers}

\blt[Overview]This section introduces notation and reviews the core concepts behind attention-based schemes in machine learning.

\subsection{Attention Heads}
\label{sec:attention}
The building blocks of transformers are attention heads. To simplify the exposition, \emph{single-head attention} is explained here but, in practice and in our experiments, an extension called \emph{multi-head attention} is used.

\begin{bullets}%
    \blt[Cross-Attention] The \emph{cross-attention operator} is,
    \begin{bullets}%
        \blt[Intuition] intuitively speaking, a function that returns a vector encoding the information that
        \blt[primitive]
        \begin{bullets}%
            \blt[Input]
            \begin{bullets}%
                \blt[Set of reference vectors]a certain set of reference vectors $\{\SidVec_1,\ldots,\SidVec_{\SidVecNum}\}\subset \rfield^{\SidVecDim}$
                \blt[Input vector]provide about a  vector $\QueVec\in \rfield^{\QueVecDim}$.
            \end{bullets}%
            \blt[Output]In particular, this operator returns a convex combination of the \emph{value} vectors $\ValFun(\SidVec_{\SidVecInd})$:
            \begin{align}
                \label{eq:valconvcomb}
                \AttFun(\SidVecMat,\QueVec) = \frac{\sum_{\SidVecInd=1}^{\SidVecNum} \AttWei
                    (
                    \SidVec_\SidVecInd, \QueVec
                    ) \ValFun(\SidVec_{\SidVecInd})}{
                    \sum_{\SidVecInd=1}^{\SidVecNum} \AttWei(
                    \SidVec_\SidVecInd, \QueVec
                    )
                }\in \rfield^{\ValVecDim},
            \end{align}
            where
            \begin{bullets}
                \blt $\SidVecMat\define [\SidVec_1,\ldots,\SidVec_{\SidVecNum}] \in \rfield^{\SidVecDim \times \SidVecNum}$
                \blt and $\AttWei(
                    \SidVec, \QueVec
                    )\geq 0$ are the so-called (unnormalized) \emph{attention weights}.
            \end{bullets}
            \blt[values] The value vectors are provided by the learnable function $\ValFun:\rfield^{\SidVecDim}\rightarrow \rfield^{\ValVecDim}$, which is normally linear. Vector $\ValFun(\SidVec_{\SidVecInd})$ encodes the information in $\SidVec_{\SidVecInd}$  that is relevant for the task at hand.

            \blt[weights] The attention weights are determined by the relation between $\QueVec$ and the vectors in $\SidVecMat$. The learnable function $\AttWei: \rfield^{\SidVecDim}\times \rfield^{\QueVecDim}\rightarrow \rfield_{++}$ can be thought of as quantifying the similarity between $\SidVec$ and $\QueVec$. In this way, if $\QueVec$ is very similar to $\SidVec_{\SidVecInd_0}$ for some $\SidVecInd_0$ and dissimilar to the remaining reference vectors, then $\AttWei
                (
                \SidVec_{\SidVecInd_0}, \QueVec
                )$ will dominate and $\AttFun(\SidVecMat,\QueVec)\approx\ValFun(\SidVec_{\SidVecInd_0})$.
            \begin{bullets}%
                \blt[inner-product attention]
                \begin{bullets}%
                    \blt[def]
                    Usually,  $\AttWei$ is the so-called  inner-product attention function\footnote{A factor of $1/\sqrt{\KeyVecDim}$ is often explicitly included inside the exponential but  it  is absorbed here into either $\KeyFun$ or $\QueFun$ to simplify notation.}
                    \begin{journalonly}
                        \begin{align}
                            \AttWei(\SidVec,\QueVec)= \exp\left[
                                \KeyFun\transpose(\SidVec) \QueFun(\QueVec)
                                \right],
                        \end{align}
                    \end{journalonly}
                    \begin{conferenceonly}
                        $\AttWei(\SidVec,\QueVec)= \exp\left[
                                \KeyFun\transpose(\SidVec) \QueFun(\QueVec)
                                \right]$,
                    \end{conferenceonly}
                    ~where $\KeyFun:\rfield^{\SidVecDim}\rightarrow \rfield^{\KeyVecDim}$ and $\QueFun:\rfield^{\QueVecDim}\rightarrow \rfield^{\KeyVecDim}$ are (typically linear) learnable functions that respectively return the so-called \emph{key} and \emph{query} vectors.
                    \blt[intuition]One can think of
                    $\KeyFun(\SidVec)$ as a vector that encodes the information in  $\SidVec$  and of  $                         \QueFun(\QueVec)$ as a vector that encodes the information relevant to $\QueVec$. Thereby, $\AttWei(\SidVec,\QueVec)$ captures how relevant $\SidVec$ is to~$\QueVec$.

                    \blt[matrix notation]By letting\footnote{Matrices in the literature on transformers  are the result of transposing the matrices here. We adopt the common notation in our community. } $\ValMat\define[\ValFun(\SidVec_{1}),\ldots,\ValFun(\SidVec_{\SidVecNum})]$ and $\KeyMat\define[\KeyFun(\SidVec_{1}),\ldots,\KeyFun(\SidVec_{\SidVecNum})]$, expression \eqref{eq:valconvcomb} becomes
                    \begin{align}
                        \label{eq:valconvcombvecmat}
                        \AttFun(\SidVecMat,\QueVec) = \ValMat \SofMax\left( \KeyMat\transpose \QueFun(\QueVec) \right).
                    \end{align}

                \end{bullets}%

            \end{bullets}

        \end{bullets}%

        \blt[matrix extension]The cross-attention operator $\AttFun$ can be extended to  matrices of input vectors $\InpMat\define[\InpVec_1,\ldots,\InpVec_{\QueVecNum}]$. In this case, $\AttFun(\SidVecMat,\InpMat)$ is a matrix whose $\QueVecInd$-th column is given by $\AttFun(\SidVecMat,\InpVec_{\QueVecInd})$ or, equivalently
        $
            \AttFun(\SidVecMat,\InpMat) = \ValMat \SofMax\left( \KeyMat\transpose \QueMat \right),$
        where $\QueMat\define[\QueFun(\QueVec_{1}),\ldots,\QueFun(\QueVec_{\QueVecNum})]$.
        \blt[self-attention]The \emph{self-attention} operator is the special case where $\SidVecMat=\InpMat$. It will be denoted as $\AttFun(\InpMat)\define \AttFun(\InpMat,\InpMat)$.

    \end{bullets}%

\end{bullets}%

\subsection{Transformers}
\newcommand{\LayNorFun}{\hc{l}}
\newcommand{\MlpFun}{\hc{f}_{\text{MLP}}}

\blt[def]Transformers~\cite{vaswani2023attention} are in essence feed-forward deep neural networks that involve attention heads.
\blt[blocks]The original architecture, here presented with some simplifications, is the composition of \emph{attention blocks}.
To introduce what an attention block is,
\begin{bullets}%
    \blt[layer normalization]define the \emph{layer normalization} operator (see references in \cite{dosovitskiy2021image}) as
    \begin{journalonly}
        \begin{align}
            \LayNorFun(\QueVec) \define a \frac{\QueVec - x_{\text{mean}}}{x_{\text{std}}} + b,
        \end{align}
    \end{journalonly}
    \begin{conferenceonly}
        $ \LayNorFun(\QueVec) \define a ({\QueVec - x_{\text{mean}}})/{x_{\text{std}}} + b$,
    \end{conferenceonly}
    ~where $a$ and $b$ are learnable parameters and $x_{\text{mean}}$ and $x_{\text{std}}$ are the sample mean and sample standard deviation of the entries of $\QueVec$. When applied to a matrix, $\LayNorFun$ operates column-wise.

    \blt[attention block]With this notation, an attention block $\AttBloFun$ is described by:
    \begin{salign}
        \label{eq:transformerblock:1}
        \InpMat'            & \define \InpMat + \AttFun(\LayNorFun_1(\InpMat))   \\
        \AttBloFun(\InpMat) & \define \InpMat' + \MlpFun(\LayNorFun_2(\InpMat')),
    \end{salign}
    where $\MlpFun$ is a multi-layer perceptron applied separately to each column, and functions $\LayNorFun_1$ and $\LayNorFun_2$ implement layer normalization without sharing their $a$ and $b$ parameters.

\end{bullets}

\section{Attention Based RME}
\label{sec:storm}

\blt[Overview] Transformers were originally proposed in the context of natural language processing (NLP). Subsequently, they were adapted to image processing~\cite{dosovitskiy2021image}. To the best of our knowledge, they have not been applied to the  interpolation of
measurements collected across space. However, this paper shows that this is not only possible but it also results in an elegant and effective radio map estimator. Sec.~\ref{sec:simulations} will show that it actually beats the state-of-the-art in RME.

\blt[functional form]A simple possibility to design a transformer-based RME estimator would be to consider the existing (grid-aware) DNN-based estimators and replace the DNN therein with a vision transformer~\cite{dosovitskiy2021image}. However, this would suffer from the limitations of these estimators; cf. Sec.~\ref{sec:intro}. Instead, an alternative route is taken here, which proceeds by adopting a somehow abstract perspective. In particular, note from the problem formulation in Sec.~\ref{sec:problem} that any estimator of  $\pow(\loc)$ given the data $\left\{(\loc_{\measind},
    \powmeas[\measind])\right\}_{\measind=1}^{\measnum}$ is a function of the form $\powest(\loc; \left\{(\loc_{\measind},
    \powmeas[\measind])\right\}_{\measind=1}^{\measnum} )$.

\blt[Straightforward approach]Note that two desirable properties for any estimator like this are (i) that it is  invariant to permutations of the measurements and (ii) that  it  accommodates an arbitrary $\measnum$. Since the cross-attention operator $\AttFun$ satisfies these two properties, one could think of constructing the feature vectors $\SidVec_{\SidVecInd}=[ \powmeas[\measind],\loc_{\measind}\transpose ]\transpose$ and estimate the map at $\loc$ as $\powest(\loc) = \AttFun(\SidVecMat,\loc)$, where $\SidVecMat\define[\SidVec_1,\ldots,\SidVec_{\measnum}]$. This could be composed with other layers to form a deep network.
\blt[missing invariances]Unfortunately, such an approach does not satisfy desirable invariance properties, as discussed next.

\subsection{Feature Design}
\label{sec:featuredesign}
\begin{bullets}%

    \blt[overview]Maxwell's equations dictate that RME exhibits certain invariance properties, such as translation and rotation invariance. Note that if the estimation of a map as a whole is considered rather than the estimation of a map at a single location $\loc$, these invariances become equivariances.
    If the invariances of a problem  are not imposed by the estimator architecture, they must be learned from data, which increases drastically the amount of  training data required to attain a target performance. This motivates a feature design that  enforces these invariances.

    \blt[Translation invariance]Translation invariance means that translating the coordinate system shall not change the estimate of $\pow(\loc)$. To impose this invariance, one can replace the feature vectors with translation-invariant features, for example $[ \powmeas_{\measind},(\loc_{\measind}-\loc)\transpose ]\transpose$. Since this would also translate the second input of $\AttFun$ to the origin ($\loc-\loc=\bm 0$), it is more appropriate to use self-attention. This means that $\powest(\loc)$ can be taken to be a function of  $\AttFun(\FeaVecMat)$, or even $\AttBloFun(\FeaVecMat)$, where $\FeaVec_{\FeaVecInd}=[ \powmeas_{\measind},(\loc_{\measind}-\loc)\transpose ]\transpose$.

    \blt[Rotation invariance]
    \begin{bullets}%
        \blt[overview] Similarly, rotation invariance means that the estimate of $\pow(\loc)$ shall not change if the coordinate system is rotated. To accommodate this invariance,
        \blt[rotate meas]the centered measurement locations $\{\loc_{\measind}-\loc\}_{\measind=1}^{\measnum}$ will be suitably rotated,
        which means that the feature vectors become $\FeaVec_{\FeaVecInd}=[ \powmeas[\measind],\RotMat(\loc_{\measind}-\loc)\transpose ]\transpose$, where $\RotMat$ is a rotation matrix.
        This rotation  is defined by one angle if  $\areadim=2$ and by two angles if $\areadim=3$. Note that rotating all locations by a certain angle amounts to rotating the coordinate system by the opposite angle. Thus, one can define a rotation by the direction in which the x-axis points after the rotation.
        \blt[rotation angle]One possibility is to choose the direction of a specific measurement location, for instance the one corresponding to the largest $\powmeas[\measind]$. However, this means that a small change in the measurements could result in a large change in the rotation angle if the index of the strongest measurement changes. Since this would render the estimator unstable, a more robust approach is adopted here, where the x-axis is rotated so that it points in the direction of
        \begin{align}
            \label{eq:rotationangledirpt}
            \sum_{\measind=1}^{\measnum} \exp( \powmeas[\measind])
            (\loc_{\measind}-\loc).
        \end{align}
        The exp  yields positive weights if $\powmeas[\measind]$ is in dB units.

    \end{bullets}%

    \blt[Additional features]Finally, besides imposing invariances, a suitable feature design can also  facilitate learning. For this reason, other features can be appended to the aforementioned feature vectors. For example, one can concatenate the cylindrical or spherical coordinates of $\RotMat(\loc_{\measind}-\loc)$ and even the sines and cosines of the resulting angular coordinates.

\end{bullets}%

\subsection{Dataset Preparation}

\blt[goal]The  data to train and test an RME estimator  typically consists of a collection of \emph{measurement sets} (MSs). Each MS contains geolocated measurements corresponding to a different \emph{true} $\Map$. For example, each MS can be collected  in a different geographical area or for different transmitter locations. MSs are then used to generate sets of training and testing \emph{examples}. Since the proposed estimator is gridless, the procedure differs from the one used in existing DNN estimators.

\blt[examples]In particular, one can proceed as follows to generate the $\TimInd$-th example. First, randomly select one of the MSs. Among the measurements in MS, choose one as the target and $\measnum\TimNot{\TimInd}$ of them as the input. As discussed later, the value of $\measnum\TimNot{\TimInd}$ must be selected depending on the training or testing goals.
With  the $\measnum\TimNot{\TimInd}$ measurements, their locations, and the location of the target measurement, construct the feature matrix $\FeaVecMat\TimNot{\TimInd}$ as indicated in Sec.~\ref{sec:featuredesign}. By denoting the target measurement as $\powmeas\TimNot{\TimInd}$, the dataset can be expressed as  $\left\{(\FeaVecMat\TimNot{\TimInd},\powmeas\TimNot{\TimInd})\right\}_{\TimInd=1}^{\TimNum}$.

\subsection{Architecture}
\label{sec:rmeest:architecture}
\blt[layers]The proposed estimator adopts a transformer architecture and comprises the blocks in the shaded area of Fig.~\ref{fig:activesensingarch}. The rest of the blocks are used for active sensing and are described later. The
feature vectors $\FeaVec_{\FeaVecInd}$ are first passed separately through a
linear layer $L$ that increases their dimension to $\EmbDim$. Then, a composition of
attention blocks is applied. The input and output of every block are vectors of
dimension $\EmbDim$. Finally, the output vectors of the last attention block are
passed through a linear layer that reduces their dimension to 1. The returned
$\measnum$ scalars can be collected in the vector
$\NetFun(\FeaVecMat)\define[\NetFun_1(\FeaVecMat),\ldots,\NetFun_\measnum(\FeaVecMat)]\transpose$.
\blt[no masking] To obtain an estimate $\powest(\loc)$,  these $\measnum$ scalars could be reduced into a single one, e.g. by averaging, and train by minimizing the mean square~error:
\begin{align}
    \frac{1}{\TimNum}\sum_{\TimInd=1}^{\TimNum} \left(\powmeas\TimNot{\TimInd}-\frac{\bm 1\transpose \NetFun(\FeaVecMat\TimNot{\TimInd})}{\measnum\TimNot{\TimInd}}\right)^2.
\end{align}

The limitation of this approach is that examples for all the necessary  values of $\measnum\TimNot{\TimInd}$ need to be included in the dataset. This may result in an unnecessarily  large dataset.
\blt[masking]To alleviate this issue,   \emph{causal self-attention} is commonly used
in the context of transformers. As discussed later, it is not fully suitable for RME, but the training time reduction may pay off.
\begin{bullets}%
    \blt[description]A causal self-attention head $\CauAttFun$ is similar to the self-attention head introduced in Sec.~\ref{sec:attention}, but the $\QueVecInd$-th output vector is only allowed to depend on the input vectors $\{\InpVec_{\FeaVecInd'}\}_{\FeaVecInd'=1}^{\QueVecInd}$. In other words, the $\QueVecInd$-th column of the matrix $\CauAttFun(\InpMat)$ is $\AttFun(\InpMat_{\QueVecInd}, \InpVec_{\QueVecInd})$, where $\InpMat_{\QueVecInd}\define[\InpVec_1,\ldots,\InpVec_{\QueVecInd}]$. An attention block that uses a causal self-attention head will be denoted as $\CauAttBloFun$.

    \blt[Benefit]Thus, if one replaces all attention heads in the aforementioned architecture with causal self-attention heads, one can train the network so that
    $\NetFun_\measind(\FeaVecMat)$ is an estimate of $\pow(\loc)$ given the first $\measind$ measurements. This can be achieved with the loss
    \begin{align}
        \frac{1}{\TimNum\measnum}\sum_{\TimInd=1}^{\TimNum}
        \sum_{\measind=1}^{\measnum}
        \left(\powmeas\TimNot{\TimInd}- \NetFun_\measind(\FeaVecMat\TimNot{\TimInd})\right)^2,
    \end{align}
    where now $\measnum\TimNot{\TimInd}=\measnum$ for all $\TimInd$.
    \blt[Limitation]The main limitation of this approach is the loss of the invariance to the permutation of the measurements. This is not a problem in NLP since the tokens in a word are ordered, but it is a problem in RME since the measurements are not ordered. However, this is the price to be paid for the  reduction in training complexity.

    \begin{journalonly}
        \blt[implementation]For future reference, note that a causal self-attention head is implemented as
        \begin{align}
            \CauAttFun(\InpMat) = \ValMat \SofMax\left( \KeyMat\transpose \QueMat - \MasMat\right),
        \end{align}
        where $\MasMat$ is a mask matrix with $-\infty$ above the main diagonal and 0 elsewhere.
    \end{journalonly}

\end{bullets}%

\section{Attention-based Active Sensing}
\label{sec:activestorm}

\blt[Idea]The estimator proposed in Sec.~\ref{sec:storm} will be referred to as STORM and will be extended next to the active sensing setting.
The idea is that,  given the measurements $\{(\loc_{\measind},
    \powmeas[\measind])\}_{\measind=1}^{\measnum}$, the candidate locations $\{\loc_{\measind}\}_{\measind=\measnum+1}^{\measnum+\CanNum}$, and a target location $\loc$, STORM must return not only an estimate of $\pow(\loc)$  but also  quantify how informative a measurement at each of these candidate locations  would be to improve this estimate.

\blt[Operational time]

\begin{bullets}%

    \blt[Features] To this end, process $\{\loc_{\measind}\}_{\measind=\measnum+1}^{\measnum+\CanNum}$ identically to the measurement locations, with the same translation and rotation as  used to obtain $\InpMat$. This yields the matrix $\CanMat$, which has one row less than $\InpMat$ since it does not contain the measurement values.
    \blt[Architecture]Decomposing $\NetFun$ into an encoder part $\NetEncFun$ and a decoder part $\NetDecFun$, the architecture of STORM is extended to obtain:
    \begin{salign}
        \NetFun(\InpMat)         & = \NetDecFun(\NetEncFun(\InpMat))                      \\
        \CanFun(\InpMat,\CanMat) & = \CanDecFun(\NetEncFun(\InpMat),\CanEncFun(\CanMat)).
    \end{salign}
    Here, $\CanFun(\InpMat,\CanMat)$ is a vector with $\CanNum$ entries, where the $\CanInd$-th entry is a scalar between 0 and 1 that quantifies how informative a measurement at $\loc_{\CanInd}$ would be to improve the estimate of $\pow(\loc)$. Its encoder and decoder parts are denoted as $\CanEncFun$ and $\CanDecFun$ and pictorially described by Fig.~\ref{fig:activesensingarch} (left).

\end{bullets}%

\blt[Training time]
\begin{bullets}%
    \blt[Features] For training, the measurements at the candidate locations are used. Thus, one can construct $\InpMat$ and $\CanMat$ as above except that $\InpMat$ now contains $\CanNum$ more columns with the candidate locations and measurements. Note that the candidate locations are not used in \eqref{eq:rotationangledirpt} to obtain the rotation angle.
    \blt[Architecture] The causal self-attention operator used in
    $\NetEncFun$ and $\NetDecFun$ is modified as follows: The first $\measnum$ output columns coincide with the $\measnum$ output columns of $\CauAttFun$ from Sec.~\ref{sec:rmeest:architecture}.
    In turn, for $\CanInd>\measnum$, the $\CanInd$-th
    output column is given by $\AttFun([\InpMat_{\measnum},\InpVec_{\CanInd}], \InpVec_{\CanInd})$. This allows one to use the $\CanInd$-th
    output column of $\NetFun(\InpMat)$ as an estimate of $\pow(\loc)$ given the $\measnum$  measurements $\{(\loc_{\measind},
        \powmeas[\measind])\}_{\measind=1}^{\measnum}$
    as well as the measurement at $\loc_{\CanInd}$, but not the measurements at other candidate locations. This is therefore the estimate in the next step of the active sensing process if the $\CanInd$-th candidate location is selected. The outputs $\NetFun_{\measnum+1}(\InpMat),\ldots,\NetFun_{\measnum+\CanNum}(\InpMat)$ will be therefore referred to as the \emph{candidate estimates}. The attention block that results from this modification will be denoted as $\ModCauAttBloFun$.
    \begin{journalonly}
        \acom{explain mask?}
    \end{journalonly}
    The overall training architecture of STORM for active sensing is presented on the right side of Fig.~\ref{fig:activesensingarch}.

    \begin{figure}
        \centering
        \includegraphics[clip, trim=2.7cm 3.45cm 11.5cm 4.0cm, width=\columnwidth]{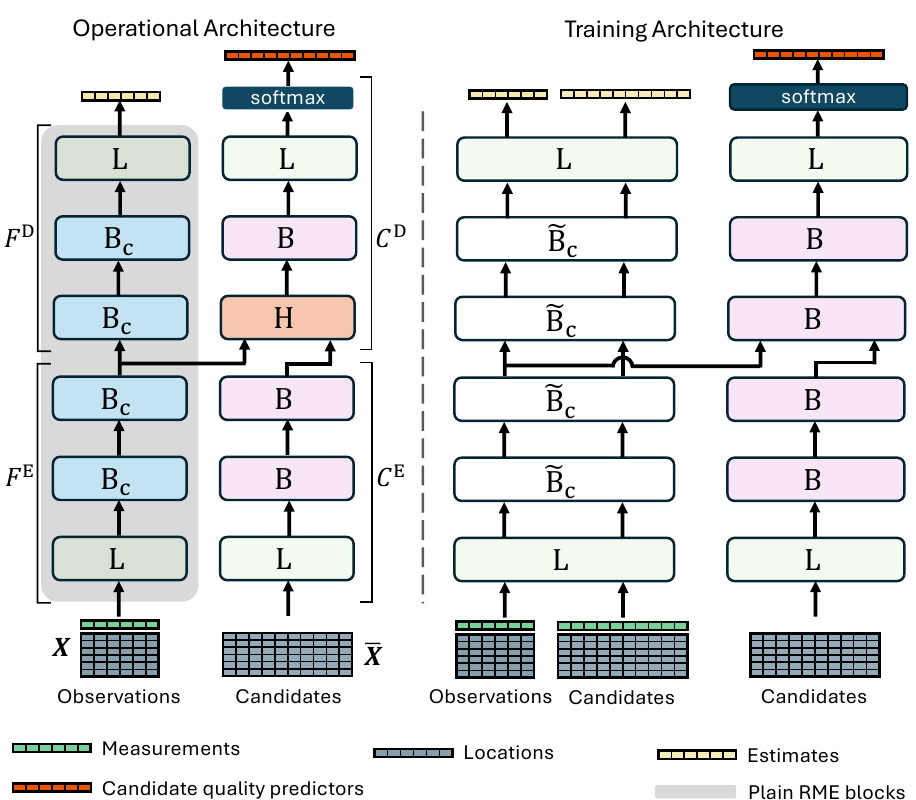}
        \caption{Architecture of STORM. $L$ denotes a linear layer, $\CauAttBloFun$ a causal self-attention block, $\ModCauAttBloFun$ a modified causal self-attention block, and $\AttFun$ a cross-attention head.}%
        \label{fig:activesensingarch}
    \end{figure}

    \blt[Loss] To train STORM  to   predict the quality of each candidate estimate without using the measurements at the candidate locations, the idea here is to construct a \emph{combined estimate} as a convex combination of the candidate estimates. The weights in this convex combination are the entries of $\CanFun(\InpMat,\CanMat)$, which do not depend on those measurements. The loss becomes:
    \begin{align}
        \frac{1}{\TimNum} & \sum_{\TimInd=1}^{\TimNum}\bigg[\frac{1}{2\measnum}
            \sum_{\measind=1}^{\measnum}
            \left(\powmeas\TimNot{\TimInd}- \NetFun_\measind(\FeaVecMat\TimNot{\TimInd})\right)^2
        \\& +\frac{1}{2}
            \left(\powmeas\TimNot{\TimInd}- \sum_{\CanInd=1}^{\CanNum}
            \CanFun_{\CanInd}(\InpMat_{\measnum}\TimNot{\TimInd},\CanMat\TimNot{\TimInd})
            \NetFun_{\measnum+\CanInd}(\FeaVecMat\TimNot{\TimInd})\right)^2
            \bigg]\nonumber.
    \end{align}
    Note that this loss also promotes good estimation performance, since it is desirable that the same network can be used both for estimation and for selecting the next measurement location,  rather than using a dedicated transformer for each task.

\end{bullets}%

    \section{Experiments with Synthetic and Real Data}
\label{sec:simulations}
This section evaluates the performance of STORM on three datasets.
\begin{bullets}%
    \blt[Compared estimators]When it comes to map estimation error, STORM is compared with three non-DNN and four DNN estimators.
    \begin{bullets}%
        \blt[Non-DNN estimators]The non-DNN estimators include K-nearest neighbors (KNN), Kriging \cite{alayafeki2008cartography,shrestha2022surveying\nextversion{,beers2004kriging,agarwal2018spectrum,boccolini2012wireless}}, and kernel ridge regression (KRR); see \cite{romero2017spectrummaps} and references therein. These estimators are trained as described in \cite{shrestha2024empiricaljpaper}.
        %
        \blt[DNN estimators]The compared DNN estimators include
        \begin{bullets}%
            \blt~[DNN 1] the completion autoencoder in \cite{teganya2020rme},
            %
            \blt~[DNN 2] the U-Net from \cite{levie2019radiounet},
            %
            \blt~[DNN 3] the U-net from \cite{krijestorac2020deeplearning},
            %
            \blt~and [DNN 4] the autoencoder in \cite{shrestha2022surveying}.
        \end{bullets}%
        \blt[proposed transformer]
        \begin{bullets}
            \blt[architecture] Unless stated otherwise, STORM uses multi-head attention with 2 heads and  embedding dimension $\ValVecDim = 48$, which
            \blt[complexity] results in around 100 k parameters.
        \end{bullets}

    \end{bullets}%
    \blt[training and testing data]Each of the considered datasets comprises multiple MSs  collected in an  $\mapSizeX\times\mapSizeY$ rectangular environment.
    \begin{bullets}%
        \blt[training testing split]These MSs are split into training and testing MSs.
        \blt[patch]To obtain each training example, one can collect the measurements inside an $\patchSize\times\patchSize$ square patch drawn uniformly at random and included in a training MS selected uniformly at random. Likewise for the testing examples.
        %
        \blt[gridless]To favor the competing (grid-aware) DNN estimators, the aforementioned square patches are aligned with the grid in which the measurements are collected. Worse performance of these estimators is expected otherwise. The grid spacing is denoted as  $\patchGridSpace$.
    \end{bullets}%

    \blt[performance metric - rmse] At each Monte Carlo iteration, the measurements
    $\left\{(\loc_{\indObs},\Mea_{\indObs})\right\}_{\indObs=1}^{\numPatchMeas}$ inside a patch are  split into two subsets by partitioning the index set
    $\patchSetMeasInd\define\{1,2,\ldots,\numPatchMeas\}$ into $\setObs$ and
    $\setNobs$, that is, $\setObs\cup\setNobs=\patchSetMeasInd$ and
    $\setObs\cap\setNobs=\emptyset$. The cardinality
    $\numPatchObs\define\left|\setObs\right|$ is fixed and presented on the horizontal axis of the figures.
    \begin{bullets}%
        \blt[input] The measurements with indices in $\setObs$ are passed to each estimator and
        \blt[output] the returned map estimate $\estimate$ is evaluated at the locations $\{\loc_{\indObs}\}_{\indObs\in\setNobs}$.
        \blt[metric] The root mean square error (RMSE) is then defined as
        \begin{align}
            \rmse\define\sqrt{\expectation\left[\frac{1}{\left|\setNobs\right|}\sum_{\indObs\in\setNobs}\left|\Mea_{\indObs}-\estimate\left(\loc_{\indObs}\right)\right|^2\right]},
        \end{align}
        where the expectation is over patches and realizations of $\setObs$.
    \end{bullets}%

    \blt[Experiments]
    \begin{bullets}%
        %

        \blt[MC: RMSE vs. num obs \ra 3 figs]
        \begin{bullets}%
            %
            \blt[ray tracing]
            \begin{bullets}%
                \blt[dataset]
                \begin{bullets}%
                    %
                    \blt The first dataset in this paper is generated by  Remcom’s Wireless InSite ray-tracing software using a 3D model of an area in downtown Rosslyn, Virginia, with $\mapSizeX
                        \approx\mapSizeY\approx 700$ m.
                    \blt Each  MS corresponds to a different  transmitter   location; see \cite{teganya2020rme} for details.
                \end{bullets}%
                %
                %
                \blt[figure] Fig.~\ref{fig:ray-tracing} shows the RMSE vs.  $\obsnum$.
                \begin{figure}
                    \centering
                    \includegraphics[width=.8\columnwidth]{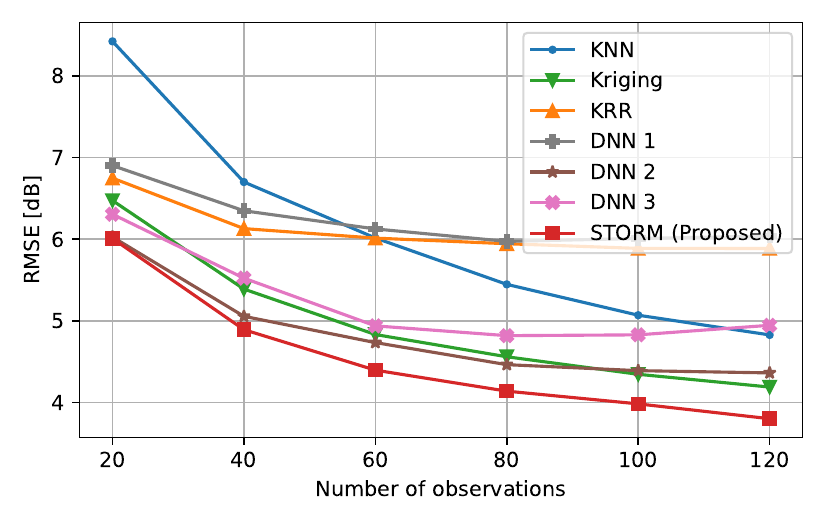}
                    \caption{ RMSE for ray-tracing data vs.  $\obsnum$   when $\patchSize= 64$ m, $\patchGridSpace=4$ m, and the estimators are trained with $\obsnum\in[20, 100]$.}
                    \label{fig:ray-tracing}
                \end{figure}
                \blt[observations] It is  observed that
                \begin{bullets}%
                    \blt STORM outperforms all other benchmarks for all $\obsnum$ despite the fact that the complexity of most competitors is significantly higher.
                    %
                \end{bullets}
            \end{bullets}%

            %
            \blt[usrp]
            \begin{bullets}%
                \blt[dataset]
                \begin{bullets}%
                    \blt The second experiment uses the USRP dataset  collected in \cite{shrestha2023empirical}.
                    \begin{journalonly}
                        The  transmitter is stationed on the ground  whereas the
                        receiver is mounted on a UAV that traverses the mapped region through a fixed trajectory. Each of the 18 MSs is generated by placing the transmitter antenna in different orientations, positions, and by introducing metallic reflectors around the transmitter.
                    \end{journalonly}
                    \blt Each MS has $\mapSizeX\approx\mapSizeY\approx53$ m and consists of approximately 12000 measurements.
                \end{bullets}%
                %
                %
                \blt[figure] Fig.~\ref{fig:usrp} shows the RMSE vs. $\obsnum$.
                \begin{figure}
                    \centering
                    \includegraphics[width=.8\columnwidth]{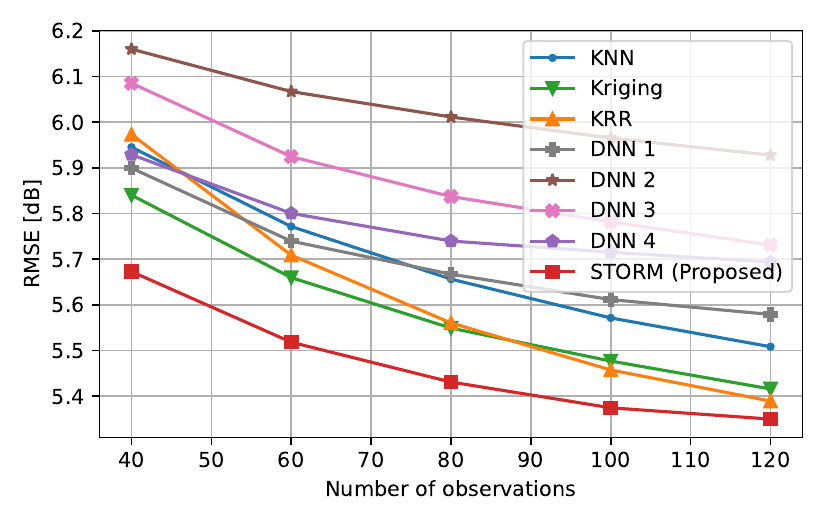}
                    \caption{ RMSE for USRP data vs. the $\obsnum$ when $\patchSize= 38.4$ m, $\patchGridSpace=1.2$ m, and the estimators are trained with $\obsnum\in[40, 100]$.}
                    \label{fig:usrp}
                \end{figure}
                \blt[observations] It can be again observed that
                \begin{bullets}
                    %
                    \blt STORM outperforms all
                    benchmarks.
                    \blt Kriging is the second best estimator, but recall that its complexity per point estimate is cubic in $\obsnum$, whereas the complexity of STORM is quadratic.
                \end{bullets}
            \end{bullets}%

            %
            \blt[4g gradiant]
            \begin{bullets}%
                \blt[dataset]
                \begin{bullets}%
                    \blt The third experiment relies on the 4G dataset from \cite{shrestha2024empiricaljpaper}.
                    \blt The transmitters are the base stations deployed by a cellular operator in a real-world 4G
                    network.
                    \begin{journalonly}
                        Conventional sector antennas are used.
                        \blt The receiver is a communication module on board a UAV.
                    \end{journalonly}
                    \begin{journalonly} with a built-in RTK  localization module, which ensures highly accurate geo-referencing of the LTE data.\end{journalonly}
                    \blt Each MS is collected in a different rectangular area with $\mapSizeX=252$ m and $\mapSizeY=260$ m.
                \end{bullets}%
                %
                %
                \blt[figure] Fig.~\ref{fig:4g-gradiant} shows the RMSE vs. $\obsnum$ for this dataset.
                \begin{figure}
                    \centering
                    \includegraphics[width=.8\columnwidth]{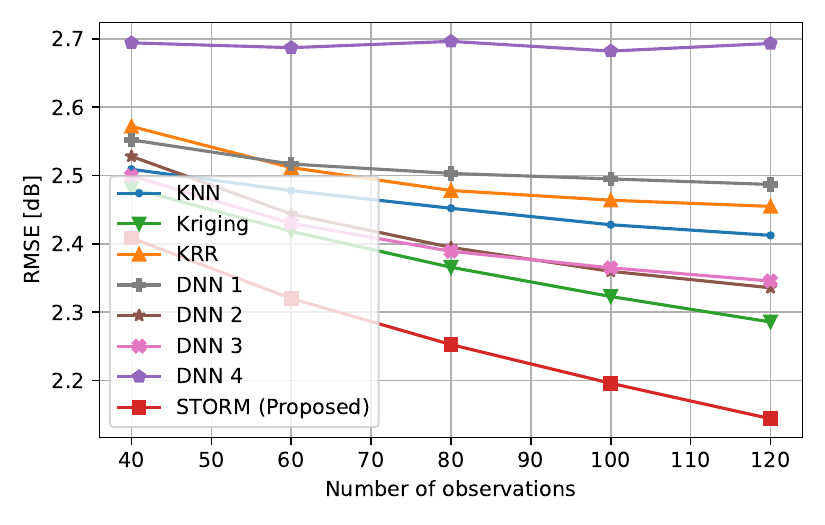}
                    \caption{ RMSE for 4G data vs. $\obsnum$ when $\patchSize= 64$ m, $\patchGridSpace=4$ m, and the estimators are trained with $\obsnum\in[20, 100]$.}
                    \label{fig:4g-gradiant}
                \end{figure}
                \blt[observations]
                \begin{bullets}%
                    \blt[outpeform]Once more, STORM offers the best estimation performance.
                    %
                    \blt[masking works]Despite the fact that STORM  is trained with $\obsnum=100$, it still
                    performs well at 120 measurements, which  corroborates the value of the considered causal attention blocks.
                    %
                \end{bullets}%
            \end{bullets}%
        \end{bullets}%

        \blt[Active Sensing: error vs. num obs \ra 1 fig] The last experiment quantifies the performance of  STORM when it comes to active sensing. For each patch, $\obsnum$ measurements, one evaluation location $\loc$, and the remaining $\CanNum=\numPatchMeas-\obsnum$ measurements are passed to STORM. The measurement with the greatest value of the quality predictor $\CanFun_{\CanInd}$ is then also given to STORM, which provides a refined estimate  of $\pow(\loc)$. This estimate is used to compute the RMSE and compared with the RMSE that results from selecting the additional measurement uniformly at random among the candidate locations. Fig.~\ref{fig:active-sensing} shows the RMSE vs. $\obsnum$. It is observed that choosing the next measurement as dictated by STORM leads to a significant improvement in the RMSE in both the ray-tracing and USRP datasets. The case of the 4G dataset is also similar but omitted due to lack of space.
        \begin{figure}
            \centering
            \includegraphics[width=.8\columnwidth]{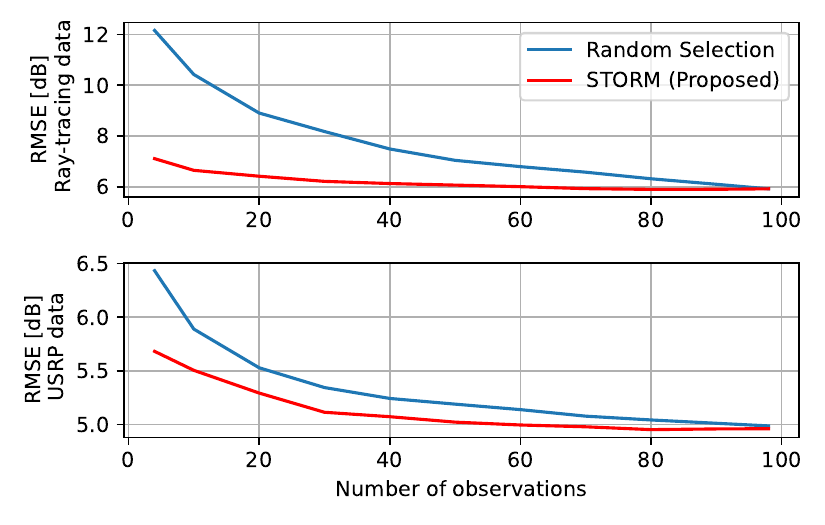}
            \caption{RMSE vs. $\obsnum$ for the active sensing problem with ray-tracing and USRP data. $\ValVecDim=20$.}
            \label{fig:active-sensing}
        \end{figure}
    \end{bullets}%
\end{bullets}%

    \section{Conclusions}
    \label{sec:conclusions}
    This paper proposed STORM, a transformer network for radio map estimation.
    This estimator operates in a gridless fashion, which circumvents many of the
    limitations of existing DNN-based estimators. It is also seen to outperform the state-of-the-art estimators in three datasets. An extension of STORM to active sensing was also proposed and seen to yield satisfactory results.

\end{bullets}%

\printmybibliography

\end{document}